# Robust Traffic Control Using a First Order Macroscopic Traffic Flow Model

Hao Liu , Christian Claudel, and Randy Machemehl

*Abstract*—Traffic control is at the core of research in transportation engineering because it is one of the best practices for reducing traffic congestion. It has been shown in recent years that the traffic control problem involving Lighthill-Whitham-Richards (LWR) model can be formulated as a Linear Programming (LP) problem given that the corresponding initial conditions and the model parameters in the fundamental diagram are fixed. However, the initial conditions can be uncertain when studying actual control problems. This paper presents a stochastic programming formulation of the boundary control problem involving chance constraints, to capture the uncertainty in the initial conditions. Different objective functions are explored using this framework, and the proposed model is validated by conducting case studies for both a single highway link and a highway network. In addition, the accuracy of relaxed optimal results is proved using Monte Carlo simulation.

## I. Introduction

TRAFFIC congestion is a global issue which is expected to become worse in the next decades as cities and populations continue to grow. As the number of vehicles increases, both people and the environment are severely affected in terms of congestion and pollution, which translates into a waste of time and money. Therefore, reducing traffic congestion is a critical issue to our human society. There are a number of ways in which congestion could be reduced, such as extending road infrastructure and decreasing user demand. However, these methods are usually expensive or sometimes impractical. Compared to these strategies, traffic flow control is a more cost-effective alternative since it aims to manage traffic operation in the best possible way given the available road facilities.

Although traffic flow modeling and controlling have attracted considerable amount of research efforts, they still remain challenging due to the random nature of many factors, such as traffic density [1], time headway [2], road capacity [3]–[5], travel demand [3], [6] and the relationship between parameters [7]. For traffic flow modeling, a common way to construct a stochastic traffic flow model is adding 'noise' into a selected deterministic conservation model and a fundamental diagram (FD), such as [8], [9]. Although such models are appropriate for uncertainties of exogenous stochastic parameters, they are incapable of dealing with endogenous uncertainties, such as driver behavior. To address this issue, Jabari and Liu [2] proposed a model to study the effect of uncertainties in time headway. Aside from the concentration on these macroscopic parameters, there are also models that research the uncertainties at a more detailed level. For example, Boel and Mihaylova [10] developed a extended cell transmission model regarding sending and receiving functions as random variables because both the location and speed of each single vehicle are random.

For robust traffic flow control, travel demand and capacity are the most commonly investigated sources of uncertainties. Zhu and Ukkusuri [3] proposed a robust speed limit control model accounting for random demand and capacity, which follow lognormal probability distributions. Then, the speed limit control problem is solved as a Markov Decision Process (MDP). Zhang and Prieur [11] proposed a two-mode Markov jump linear hyperbolic (MJLH) model and designed an on-ramp metering control method for freeway traffic governed by AW-Rascle traffic flow model [12] and Greenshields FD [13]. In this model, the transition of two traffic modes, i.e., free flow mode and congestion model, is modeled as a Morkov process. Besides freeway networks, uncertainties also play an important role for signal timing designs including fixed time control [14] and adaptive signal control [15], [16]. All these models focus on the uncertainties in demand and assume inflow rates are random variables following known probability distributions.

Traffic densities are another key parameter for traffic flow control methods, but there do not exist practical or inexpensive ways to measure this directly on the field. In order to estimate the densities, indirect measurements, such as speed and volume counts, and an appropriate conversion model are usually resorted to. As a result, there may exist uncertainties in the density estimation due to random errors of sensors and systematic errors from the conversion model. However, unlike other random parameters introduced in previous paragraphs, little attention has been attracted by the uncertainties in initial densities. Although both deterministic approaches [17], [18] and stochastic approaches [19]–[22] have been proposed to make the density estimation more dependable, their goal is to reconstruct the evolution of traffic states with a high accuracy rather than study the impact of uncertainties. To fill this gap, a robust control model is proposed in this paper to deal with this type of randomness.

Traffic flow is usually modeled by deterministic Partial Differential Equations (PDEs) [12], [17], [23], [24]. The Lighthill-Whitham-Richards (LWR) model [23], [24], which is the most commonly used macroscopic traffic flow model, and a triangular FD are used in this paper to illustrate the traffic dynamics. The proposed robust control model is based on the Lax-Hopf solution, derived by Mazaré *et al.* [25], to the LWR model. This solution outperforms other solutions in two ways. First, contrary to the models in [26]–[29] that discretize PDEs into ordinary differential equations (ODEs), it is an analytical solution to the LWR PDE and does not require any discretization or approximation. Second, this solution is grid-free which can calculate the traffic state at any point directly from the initial and boundary conditions without any knowledge of prior event. Other analytical solutions such as front tracking method [30] require full knowledge of prior events and may have exponential growth of waves in some situations as waves "bounce" back and forth from the boundary conditions, which results in low efficiency.

Built upon the Lax-Hopf solution, this paper proposes a stochastic traffic control model to handle the uncertainties in initial densities by using chance constraints, which is a major approach to solving optimization problems with random parameters, see [6]. Usually, as long as the real densities are higher than the value used in a classical control model ignoring the randomness, queue spillover is inevitable. Contrarily, the chance constraints ensure robust optimal control is feasible with a high probability. The feasibility in this framework means the road links have enough space to accept the inflows derived from the model. In other words, the proposed model can avoid queue spillover for a wide range of realization of initial densities.

The rest of this paper is organized as follows. Section II reviews the Lax-Hopf solution to the LWR model and presents the compatibility conditions need to satisfied by the solution. Section III introduces the uncertainties of initial densities to the compatibility conditions through chance constraints. Then, the chance constraints are linearized, and a case study for a single highway link is conducted. Section IV exhibits a case study for a highway network for validation. It shows that unlike the non-robust control model, the proposed model is able to avoid congestion and maintain the network mobility with random initial densities. The accuracy of the relaxed solutions are demonstrated through Monte Carlo simulations in Section V. Section VI summarizes the work and discusses potential future research.

## II. MOSKOWITZ SOLUTIONS FROM LAX-HOPF FORMULA

This section covers the Lax-Hopf solutions to the LWR model and the constraints need to be satisfied. The proposed model is built on this framework. The organization of this section is as follows: part A introduces the traffic flow model; by using Lax-Hopf formula, part B shows the Moskowitz solutions, derived in [25], [31], associated with each value condition; such solutions have to satisfy the compatibility condition to ensure that the true solutions equal the true value conditions at corresponding points, part C presents the mathematical form of this compatibility condition. For details regarding the derivation and proof, we refer readers to references [32]–[34].

### A. Traffic Flow Models

Lighthill-Whitham-Richards (LWR) PDE model [23], [24] is one of the most commonly used models to depict the evolution of traffic flow,

$$\frac{\partial \rho(t,x)}{\partial t} + \frac{\partial \psi(\rho(t,x))}{\partial x} = 0 \quad (1)$$

where $\rho(t,x)$ is the density of the point $x$ away from a reference point at time $t$, $\psi$ is the concave Hamiltonian, which is used to denote the experimental relationship between flow and density. For simplicity, a triangular FD is used to present the relationship between flow and density,

$$\psi(\rho) = \begin{cases} v_f \rho & \rho \in [0, \rho_c] \\ w(\rho - \rho_m) & \rho \in [\rho_c, \rho_m] \end{cases} \quad (2)$$

where $v_f$ is the free flow speed, $w$ is the congestion speed, $\rho_c$ is the critical density where the flow is maximum, $\rho_m$ is the jam density, where the flow is zero due to the total congestion. Those parameters are dependent and the relationship between them can be expressed as,

$$\rho_c = \frac{-w\rho_m}{v_f - w} \quad (3)$$

Alternatively, the traffic flow can be modeled by a scalar function $M(t,x)$, known as the Moskowitz function [35], which represents the index of the vehicle at $(t, x)$. The relationship between the Moskowitz function and density and flow can be expressed as,

$$\rho(t,x) = -\frac{\partial M}{\partial x}, \quad q(t,x) = \frac{\partial M}{\partial t} \quad (4)$$

Therefore, another traffic flow model, Hamilton-Jacobi (H-J) PDE, can be obtained from the integration of the LWR PDE model (1) in space,

$$\frac{\partial M(t,x)}{\partial t} - \psi(-\frac{\partial M(t,x)}{\partial x}) = 0 \quad (5)$$

For the purposes of this work, the spatial domain $[\xi, \chi]$, where $\xi$ is the upstream boundary and $\chi$ is the downstream boundary, were divided evenly into $k_{max}$ segments; the time domain $[0, t_{max}]$, where $t_{max}$ is the simulation time, were divided evenly into $n_{max}$ segments, as shown in Figure 1. Also, we defined $K = \{1, \ldots, k_{max}\}$ and $N = \{1, \ldots, n_{max}\}$. The piecewise affine initial condition $M_k(t,x)$, upstream boundary condition $\gamma_n(t,x)$, and downstream boundary condition $\beta_n(t,x)$ in terms of the Moskowitz function are defined as follows,

$$M_k(t,x) = \begin{cases} -\sum_{i=1}^{k-1} \rho(i)X \\ -\rho(k)(x - (k-1)X), & \text{if } t = 0 \\ & \text{and } x \in [(k-1)X, kX] \\ +\infty, & \text{otherwise} \end{cases} \quad (6)$$

$$\gamma_n(t, x) = \begin{cases} \sum_{i=1}^{n-1} q_{in}(i)T \\ +q_{in}(n)(t-(n-1)T), & \text{if } x = \xi \\ & \text{and } t \in [(n-1)T, nT] \\ +\infty, & \text{otherwise} \end{cases} \quad (7)$$

$$\beta_n(t, x) = \begin{cases} \sum_{i=1}^{n-1} q_{out}(i)T \\ +q_{out}(n)(t-(n-1)T) \\ -\sum_{k=1}^{k_{max}} \rho(k)X, & \text{if } x = \chi \\ & \text{and } t \in [(n-1)T, nT] \\ +\infty, & \text{otherwise} \end{cases} \quad (8)$$

where the subscripts $k$ and $n$ indicate the index of spatial segment and time step, respectively, $X$ and $T$ are the length of the spatial segment and time segment, respectively, $\rho(k)$ is the initial density for the $k$th spatial segment, $q_{in}(n)$ and $q_{out}(n)$ are the inflow and outflow, respectively, for the $n$th time segment at boundaries.

Since this method discretizes the initial conditions and boundary conditions and assumes Equations (6)-(8) are piecewise linear functions, both spatial segment $X$ and time step $T$ have significant impact on the results. The spatial segment should be the shortest road segment within which the initial density can be regarded as a constant. To ensure these two parameters are consistent with the characteristic velocities of traffic (on the order of $v_f$), we choose the segment length and time step $T$ such as $|v_f T/X| < 1$. Note that this condition is similar to the classical Courant-Friedrichs-Lewy (CFL) condition used for solving discretized PDEs using some first order numerical schemes, though our numerical scheme is unconditionally stable (and exact) and larger time steps are allowable.

*B. Moskowitz Solutions*

The Barron-Jensen/Frankowska (B-J/F) solution [36], [37] was incorporated in to solve the H-J equation. The B-J/F solutions are fully characterized by the Lax-Hopf formula.

*Definition 1 (Value Condition):* A value condition $c(\cdot, \cdot)$ is a lower semicontinuous function defined on a subset of $[0, t_{max}] \times [\xi, \chi]$.

In the following, all of the initial conditions and boundary conditions are regarded as value conditions.

*Proposition 1 (Lax-Hopf Formula):* Let $\psi(\cdot)$ be a concave and continuous Hamiltonian, and let $c(\cdot, \cdot)$ be a value condition. The B-J/F solution $M_c(\cdot, \cdot)$ to (5) associated with $c(\cdot, \cdot)$ is defined [38]–[40] by

$$M_c(t, x) = \inf_{(u,T) \in (\varphi^*) \times R_+} (c(t-T, x+Tu) + T\varphi^*(u)) \quad (9)$$

where $\varphi^*(\cdot)$ is the Legendre-Fenchel transform of an upper semicontinuous Hamiltonian $\psi(\cdot)$, which is given by,

$$\varphi^*(u) := \sup_{p \in Dom(\psi)} [p \cdot u + \psi(p)] \quad (10)$$

Using this Lax-Hopf formula, for any point $(t, x)$, each value condition defined in (6)-(8) generates a Moskowitz solution, as shown in Figure 1, and the corresponding explicit

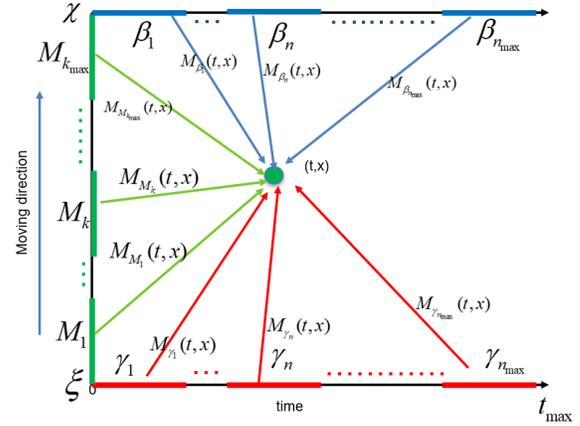

Fig. 1. Schematic of Moskowitz solution. Initial, upstream and downstream conditions shown by green, red and blue bars, respectively.

Moskowitz solutions can be expressed as (11)-(13), as shown at the bottom of the next page. Note that the Lax-Hopf formula allows for the use of any concave and piecewise linear FD, and relaxing the triangular FD assumption does not affect the formulation of optimization models in following sections.

*C. Linear Constraints*

As shown in Figure 1, each value condition generates one solution at a certain point including the domain of value conditions. Based on the Inf-morphism property, defined as Proposition 2, the real Moskowitz solution is equal to the minimum among all those solutions [38].

*Proposition 2 (Inf-Morphism Property):* Let the value condition $c(\cdot, \cdot)$ be minimum of a finite number of lower semicontinuous functions:

$$\forall (t, x) \in [0, t_{max}] \times [\xi, \chi], \quad c(t, x) := \min_{j \in J} c_j(t, x) \quad (14)$$

The corresponding solution $M_c(\cdot, \cdot)$ can be decomposed [38], [39] as

$$\forall (t, x) \in [0, t_{max}] \times [\xi, \chi], \quad M_c(t, x) := \min_{j \in J} M_{c_j}(t, x) \quad (15)$$

According to the Inf-morphism property, for any point in the domain of a value condition $c(t, x)$, the Moskowitz solutions derived from other value conditions have to be equal to or larger than $c(t, x)$ so that the real solution is equal to the value condition. This leads to the compatibility conditions defined below,

*Proposition 3 (Compatibility Conditions):* Use the value condition $c(t, x)$ and the corresponding solution in *Proposition 2*. The equality $\forall (t, x) \in Dom(c), M_c(t, x) = c(t, x)$ is valid if and only if the inequalities below are satisfied,

$$M_{c_j}(t, x) \geq c_i(t, x), \quad \forall (t, x) \in Dom(c_i), \quad \forall (i, j) \in J^2 \quad (16)$$

By substituting the Moskowitz solutions (11)-(13) and value conditions (6)-(8) into the compatibility conditions, the following constraints are obtained. For details regarding the

$$M_{M_k}(t, x) = \begin{cases} +\infty, & \text{if } x \leq (k-1)X + tw \\ & \text{or } x \geq kX + v_f t \quad (11\text{a}) \\ -\sum_{i=1}^{k-1} \rho(i)X + \rho(k)( & \text{if } x \geq (k-1)X + v_f t \quad (11\text{b}) \\ tv_f + (k-1)X - x), & \text{and } x \leq kX + v_f t \\ & \text{and } \rho(k) \leq \rho_c \\ -\sum_{i=1}^{k-1} \rho(i)X + \rho_c( & \text{if } x \leq (k-1)X + v_f t \quad (11\text{c}) \\ tv_f + (k-1)X - x), & \text{and } x \geq (k-1)X + tw \\ & \text{and } \rho(k) \leq \rho_c \\ -\sum_{i=1}^{k-1} \rho(i)X + \rho(k)( & \text{if } x \leq kX + tw \quad (11\text{d}) \\ tw + (k-1)X - x) & \text{and } x \geq (k-1)X + tw \\ -\rho_m tw, & \text{and } \rho(k) \geq \rho_c \\ -\sum_{i=1}^{k} \rho(i)X & \text{if } x \leq kX + tv_f \quad (11\text{e}) \\ +\rho_c(tw + kX - x) & \text{and } x \geq kX + tw \\ -\rho_m tw, & \text{and } \rho(k) \geq \rho_c \end{cases}$$

$$M_{\gamma_n}(t, x) = \begin{cases} +\infty, & \text{if } t \leq (n-1)T + \dfrac{x-\xi}{v_f} \quad (12\text{a}) \\ \sum_{i=1}^{n-1} q_{in}(i)T + q_{in}(n)( & \text{if } t \geq (n-1)T + \dfrac{x-\xi}{v_f} \quad (12\text{b}) \\ t - \dfrac{x-\xi}{v_f} - (n-1)T), & \text{and } t \leq nT + \dfrac{x-\xi}{v_f} \\ \sum_{i=1}^{n} q_{in}(i)T + \rho_c v_f( & \text{otherwise} \quad (12\text{c}) \\ t - \dfrac{x-\xi}{v_f} - nT), & \end{cases}$$

$$M_{\beta_n}(t, x) = \begin{cases} +\infty, & \text{if } t \leq (n-1)T + \dfrac{x-\chi}{w} \quad (13\text{a}) \\ -\sum_{k=1}^{k_{max}} \rho(k)X + & \text{if } t \geq (n-1)T + \dfrac{x-\chi}{w} \quad (13\text{b}) \\ \sum_{i=1}^{n-1} q_{out}(i)T + & \text{and } t \leq nT + \dfrac{x-\chi}{w} \\ q_{out}(n)(t - \dfrac{x-\chi}{w} & \\ -(n-1)T) - & \\ \rho_m(x-\chi), & \\ -\sum_{k=1}^{k_{max}} \rho(k)X & \text{otherwise} \quad (13\text{c}) \\ +\sum_{i=1}^{n} q_{out}(i)T + & \\ \rho_c v_f(t - nT - \dfrac{x-\chi}{v_f}), & \end{cases}$$

derivation and proof, the readers are referred to [25], [31].

$$\begin{cases} M_{M_k}(0, x_p) \geq M_p(0, x_p) & \forall (k, p) \in K^2 \\ M_{M_k}(pT, \chi) \geq \beta_p(pT, \chi) & \forall k \in K, \quad \forall p \in N \\ M_{M_k}(\frac{\chi - x_k}{v_f}, \chi) \geq \beta_p(\frac{\chi - x_k}{v_f}, \chi) & \forall k \in K, \quad \forall p \in N \\ \qquad \text{s.t.} \quad \frac{\chi - x_k}{v_f} \in [(p-1)T, pT] \\ M_{M_k}(pT, \xi) \geq \gamma_p(pT, \xi) & \forall k \in K, \quad \forall p \in N \\ M_{M_k}(\frac{\xi - x_{k-1}}{w}, \xi) \geq \gamma_p(\frac{\xi - x_{k-1}}{w}, \xi) & \forall k \in K, \quad \forall p \in N \\ \qquad \text{s.t.} \quad \frac{\xi - x_{k-1}}{w} \in [(p-1)T, pT] \end{cases}$$
(17)

$$\begin{cases} M_{\gamma_n}(pT, \xi) \geq \gamma_p(pT, \xi) & \forall (n, p) \in N^2 \\ M_{\gamma_n}(pT, \chi) \geq \beta_p(pT, \chi) & \forall (n, p) \in N^2 \\ M_{\gamma_n}(nT + \frac{\chi - \xi}{v_f}, \chi) \geq \beta_p(nT + \frac{\chi - \xi}{v_f}, \chi) & \forall (n, p) \in N^2 \\ \qquad \text{s.t.} \quad nT + \frac{\chi - \xi}{v_f} \in [(p-1)T, pT] \end{cases}$$
(18)

$$\begin{cases} M_{\beta_n}(pT, \xi) \geq \gamma_p(pT, \xi) & \forall (n, p) \in N^2 \\ M_{\beta_n}(nT + \frac{\xi - \chi}{w}, \xi) \geq \gamma_p(nT + \frac{\xi - \chi}{w}, \xi) & \forall (n, p) \in N^2 \\ \qquad \text{s.t.} \quad nT + \frac{\xi - \chi}{w} \in [(p-1)T, pT] \\ M_{\beta_n}(pT, \chi) \geq \beta_p(pT, \chi) & \forall (n, p) \in N^2 \end{cases}$$
(19)

These constraints add upper bounds for each value condition to make them physically feasible. For example, assume the road segment is completely congested, i.e., $\rho(k) = \rho_m, \quad \forall k \in K$, the left hand sides of the fourth and fifth constraints in (17) are computed from (11) and equal to zero for small $p$'s, which are the indices of time steps. As a result, the right hand sides have to be equal to zero. Consequently, based on Equation (12), the inflows at these steps must equal to zero which indicates that the road cannot accommodate any arrivals before some vehicles in the first segment depart.

Above all, a traffic control model derived based on the Lax-Hopf solutions needs to satisfy the compatibility conditions, see examples in [32]–[34]. In such formulations, the boundary conditions are the decision variable, i.e. the objective function can be realized through controlling the inflow and outflow on a traffic link; the objective function is a function of the decision variables such as minimization of delay and maximization of throughput. For a freeway link with an on-ramp, the proposed control method can be realized by on-ramp signals. For a general freeway link without an on-ramp, some control strategies, such as dynamic tolling, have the potential to control the flows at the boundaries of the link. Although the strategies used to control the boundary flow on such a highway link are not very mature so far, it is reasonable to assume that most highway links could be controlled in the future when connected and autonomous vehicles become prevalent. The parameters and control variables of the models proposed in this paper are summarized in Table I. All of the parameters except for initial densities are constant. Also, although outflows are also control variables in the model, they are only constrained by the compatibility conditions. In other words, they are determined by model parameters and optimal inflows, and do not need exogenous control in reality.

TABLE I
NOTATIONS

| Notation | Definition |
|---|---|
| **Parameters** | |
| ***Flow Model:*** | |
| $\psi$ | fundamental diagram |
| $v_f(j)$ | free flow speed of link $j$ (m/s) |
| $\rho_c(j)$ | critical density of link $j$ (veh/m) |
| $w(j)$ | backwave speed of link $j$ (m/s) |
| $\rho_m(j)$ | jam density of link $j$ (veh/m) |
| $c(j)$ | capacity of link $j$ (veh/s) |
| $\bar{\rho}(j)$ | mean of initial density on link $j$ |
| $\rho_s(j)$ | standard deviation of initial density on link $j$ |
| ***Link and Node:*** | |
| $n_{lane}(j)$ | number of lanes of link $j$ |
| $p(i, j)$ | turning ratio from link $j$ to link $i$ |
| ***Discretization:*** | |
| $X$ | length of spatial segment (m) |
| $T$ | time step size (s) |
| $k_{\max}(j)$ | number of spatial segments on link $j$ |
| $n_{\max}$ | number of time steps |
| **Control variables** | |
| $q_{\text{in}}(i, j)$ | inflow of link $j$ at time $i$ (veh/s) |
| $q_{\text{out}}(i, j)$ | outflow of link $j$ at time $i$ (veh/s) |
| $q_{\text{on}}(i, j)$ | inflow of on-ramp $j$ at time $i$ (veh/s) |
| $q_{\text{off}}(i, j)$ | outflow of off-ramp $j$ at time $i$ (veh/s) |

## III. ROBUST CONTROL FOR A SINGLE HIGHWAY LINK WITH UNCERTAINTY IN INITIAL CONDITIONS

The initial conditions are known and fixed in the compatibility conditions shown in the previous section. In reality, however, uncertainties exist in the initial conditions due to errors in the measurement. To deal with this situation, a stochastic programming model is derived in which the initial conditions are random variables with normal distributions in this section. Traffic state estimation has drawn much attention in the past decades due to its contribution to reducing travel delay and ensuring travel time reliability. The normal distribution is widely used to model the uncertainties in different traffic states, such as velocity field in a CTM-v model [41], cumulative vehicle counts at the ends of a highway link [42] and traffic flow [9]. Without loss of generality, it is reasonable to approximate uncertainties in initial densities as normal distributions. However, this assumption is not necessary for our model. In fact, benefiting from the monotonicity of the Moskowitz solution, our model is applicable for all general distributions as long as the cumulative density functions can be obtained.

### A. The Stochastic Programming Formula

Motivated by boundary control problems, in the rest of this paper, the objective functions are only functions of boundary conditions and there are no uncertainties in the objective functions. A general inequality form of an LP problem is

$$\min_{\substack{q_{\text{in}}, q_{\text{out}}; \\ q_{\text{on}}, q_{\text{off}}}} f(x)$$
$$\text{s.t.} \quad Ax \geq b \qquad (20)$$

where $f(x)$ is the objective function; $Ax \geq b$ indicate the constraints including compatibility conditions and other possible conditions, such as flow transition at nodes. When there is uncertainty in the constraints, we can formulate the counterpart stochastic programming problem with chance constraints, as follows:

$$\min_{\substack{q_{\text{in}}, q_{\text{out}};\\ q_{\text{on}}, q_{\text{off}}}} f(x) \\ s.t. \quad Pr\{Ax \geq b\} \geq 1 - \alpha \quad (21)$$

where $1 - \alpha$ is the confidence level of the chance constraint. Assume $\rho(k)$ is subjected to a normal distribution with mean and standard deviation of $(\rho_k, \sigma_k)$. Then, we can convert the constraints (17)-(19) to chance constraints. For example, the constraint

$$M_{M_k}(pT, \xi) \geq \gamma_p(pT, \xi), \quad \forall k \in K, \quad \forall p \in N \quad (22)$$

should be converted to,

$$P(M_{M_k}(pT, \xi) \geq \gamma_p(pT, \xi)) \geq 1 - \alpha, \quad \forall k \in K, \quad \forall p \in N \quad (23)$$

TABLE II
STATISTICS OF TRAFFIC STATES

| stations | $\bar{q}$ (vph) | $\sigma(q)$ (vph) | $\bar{v}$ (m/s) | $\sigma(v)$ (m/s) | $\rho_k$ (vpm) |
|---|---|---|---|---|---|
| 400536 | 6663.9 | 264.1 | 28.4 | 0.3 | 0.065 |
| 400488 | 4939.0 | 191.7 | 28.9 | 0.3 | 0.047 |
| 401561 | 5120.5 | 325.4 | 27.6 | 0.6 | 0.052 |
| 400611 | 5382.3 | 198.6 | 26.1 | 1.8 | 0.057 |
| 400928 | 5079.6 | 368.8 | 27.5 | 0.8 | 0.051 |
| 400284 | 5258.4 | 266.5 | 25.9 | 1.4 | 0.056 |

Thanks to the monotonicity and the piecewise linearity of the Moskowitz solutions, the chance constraints can be converted to linear constraints (38)-(42). The derivation of the corresponding linear constraints is shown in Appendix A.

It should be noticed that the FD is empirical, and the model could be more appropriate if those variables were regarded as random variables as well. In this paper, however, we only consider the uncertainties in initial conditions because of the complexity of dealing with stochastic model parameters, which would result in a non-tractable control problem. Although the Moskowitz solutions (11)-(13) are piece wise linear functions in the initial and boundary conditions, the bilinear terms of parameters in these solutions, e.g. $\rho_c t v_f$ in the solution (11c), make it hard to solve the traffic control problem when the uncertainty of parameters is introduced.

### B. Case Study for a Single Highway Link

We implemented our framework onto a single highway link with 4 lanes located between the Freeway Performance Measurement System (PeMS) vehicle detection stations 400536 and 400284 on Highway I-880 N around Hayward, CA, USA. We divided this spatial domain of 3.858 km into 6 even segments and created a temporal domain of 7 min with 21 even segments. The model parameters were defined as follows: the capacity $C = 8000$ vph; the critical density $\rho_c = 0.074$ /m; the free flow speed $v_f = 30$ m/s; the jam density $\rho_m = 0.5$ /m. Although the densities at a specific time are not directly measured, they can be inferred from the measurement of the flow and speed. Table II shows the data measured from the weekdays between 05/01/2018 and 05/31/2018, and the time interval is from 9:00 am to 10:00 am. In Table II, $\bar{q}$ and $\sigma(q)$ denote the mean value and standard deviation of the flow, respectively; $\bar{v}$ and $\sigma(v)$ represent the mean value and standard deviation of the speed, respectively. The means of density, $\rho_k$, are obtained by the approximation $\bar{q}/\bar{v}$. For this case study, we assumed the initial densities are normally distributed, and the mean values are equal to $\rho_k$.

We used the IBMIlogCplex solver in Matlab to solve the LPs. The International System of Units was adopted and the units were omitted in the following analysis for simplicity.

Four scenarios with different standard deviations (0.003, 0.006, 0.009 and 0.012) in the initial condition segments were investigated. For each single scenario, the standard deviation in all initial condition segments were the same and denoted by $\sigma$, and the confidence level, $1 - \alpha$, was equal to 97.5%. The optimization model is,

$$\max_{q_{\text{in}}, q_{\text{out}}} h \sum_{i=1}^{n_{max}} q_{out}(i) - \sum_{i=2}^{n_{max}} |q_{out}(i) - q_{out}(i-1)| \\ s.t. \quad A_{model} x \geq b_{model} \\ q_{out}(i) \geq 0 \quad \forall i \in N \quad (24)$$

where the first term in the objective function is to maximize the total outflow during the simulation time, and the second term smooths the outflows through forcing the difference between two adjacent outflows to be as small as possible. $h > 2$ ensures that the optimally total outflow is not impacted, and the proof is shown in Appendix B. $h = 3$ is employed in the rest of this paper. $A_{model}$ and $b_{model}$ are the coefficient matrix and right-hand side vector coming from the inequalities (38) to (42). The second term in the objective function was not a linear function of decision variables. To linearize the objective function, we added another variable vector $q_d(i)$, $i \in \{2, 3, \ldots n_{max}\}$ and extra constrains into this model,

$$\max_{q_{\text{in}}, q_{\text{out}}, q_d} h \sum_{i=1}^{n_{max}} q_{out}(i) - \sum_{i=2}^{n_{max}} q_d(i) \\ s.t. \quad A_{model} x \geq b_{model} \\ q_d(i) \geq q_{out}(i) - q_{out}(i-1), \\ \forall i \in \{2, 3, \ldots, n_{max}\} \\ q_d(i) \geq q_{out}(i-1) - q_{out}(i), \\ \forall i \in \{2, 3, \ldots, n_{max}\} \\ q_{out}(i) \geq 0 \quad \forall i \in N \quad (25)$$

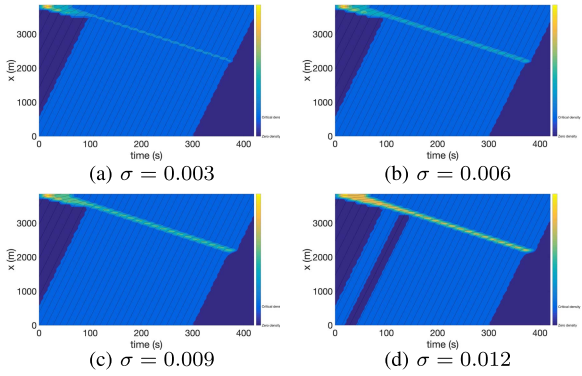

Fig. 2. Solution to robust control problem (24).

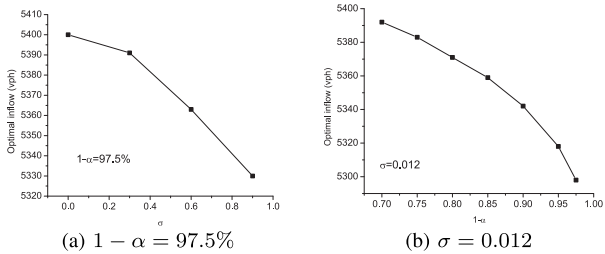

Fig. 3. Relationship between optimal total outflow and uncertainty.

The Moskowitz functions corresponding to the optimal solutions are shown in Fig. 2. All of the shockwaves, which are the congested section indicated by the yellow band, are consecutive, which makes the outflow smooth without losing optimality in terms of total outflow. In these cases, the confidence level is fixed, so the confidence interval for the initial condition is wider when the standard deviation is larger. Intuitively, the chance constraints forced the solution to satisfy (17)-(19) for all of the values of the initial conditions in the confidence interval. Therefore, the wider the confidence interval is, the more restricted the feasible region is. As this is a maximization problem, the optimal value should be lower for the case with larger standard deviation (i.e. with smaller feasible region). With increasing standard deviation, the temporal width of the shock wave (the yellow band) was wider because fewer vehicles could proceed, as shown in Fig. 2.

Note that the chance constraints are used to ensure the compatibility conditions hold with a high probability, but this does not necessarily lead to a better overall performance in terms of improving throughput. The optimal control becomes more conservative with the increase of the confidence level. Therefore, if the real initial densities are relative low, the robust optimal control may not make the best use of road capacity. To explore the influence of the standard deviation of the initial conditions and the confidence level on the optimal control, different scenarios with different variance and confidence levels were solved. The average optimal inflows are shown in Fig. 3.

In Fig. 3(a), the confidence level is fixed, and the optimal average inflow decreases as the standard deviation of the initial conditions increases. In Fig. 3(b), the standard deviation is fixed, and the optimal average inflow decreases as the confidence level increases. The trends of these two curves can be explained by the same reason as before: as variability increases, optimal inflow decreases. For the case with a confidence level of 97.5%, when the standard deviation goes up to 0.07, the LP becomes infeasible. This can be seen through Fig. 11. A large standard deviation may lead to a large $z_{1-\alpha}$, and the associated $f_2(\rho_k + z_{1-\alpha}\sigma_k)$ may be below zero because it is a decreasing function of $\rho_k$. This will generate a constraint forcing the boundary conditions to be less than zero, which will result in infeasibility. Above all, the feasible region shrinks as the variation increases.

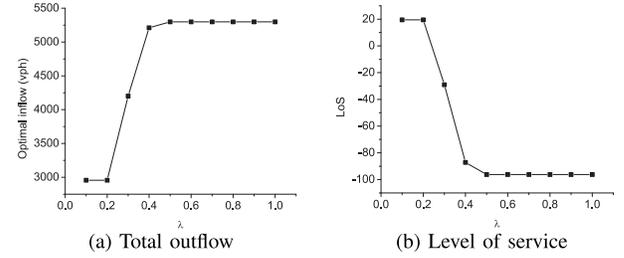

Fig. 4. Components of optimal value.

In reality, we may not only want to maximize the outflow, but also to minimize the congestion. There are several ways to realize this objective, such as adding other constraints to represent the worst level of service and change the objective function. Here, we formulated this problem as follows:

$$\min_{q_{\text{in}}, q_{\text{out}}, Q} -\lambda \sum_{i=1}^{n_{max}} q_{out}(i) + (1-\lambda)Q$$

$$s.t. \quad Q \geq \sum_{j=1}^{i}(q_{in}(j) - q_{out}(j)), \quad \forall i \in N$$

$$A_{model} x \geq b_{model}$$

$$q_{out}(i) \geq 0 \quad \forall i \in N \quad (26)$$

where $QT + \sum_{k=1}^{k_{max}} \rho(k)X$ is the maximum number of vehicles stuck in the link during the simulation, $\lambda$ and $1 - \lambda$ are the weights of total outflow and $Q$, respectively. The sum of weighted negative total outflow and $Q$ is the new objective function.

The standard deviation of the initial conditions is 0.012, and the confidence level of the chance constraints is 97.5%. To make the result more intuitive, we defined the level of service as $LoS = -QT$. Optimal inflows for different weightings (Fig. 4) show that there is a trade-off between the optimal inflows and the level of service. With the increase in $\lambda$, more vehicles can go through the highway link with a poorer level of service. When $\lambda \geq 0.5$, the sum of inflows becomes dominant, and the optimal values are stable. This critical value of $\lambda$ depends on the standard deviation of the densities and the confidence level.

### C. Computational Complexity

Without the loss of generality, let us assume inflows and outflows are the control variables, e.g., Equation (24). The number of control variables in our approach is equal to $2 \times n_{\max}$ no matter the number of spatial segments on the link.

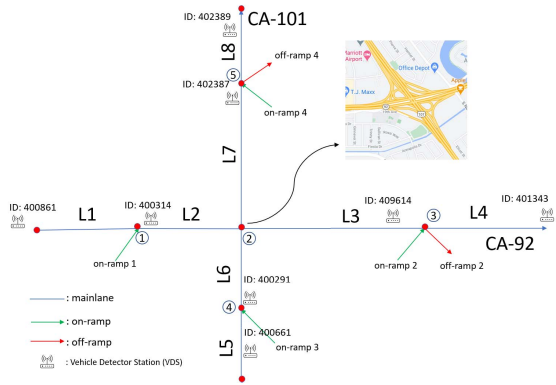

Fig. 5. Highway network layout.

To the best of our knowledge, many optimization models that discrete links into cells require at least $n_{max} \times k_{max}$ decision variables due to the traffic propagation between cells [43]. Therefore, the proposed approach is expected to outperform such methods whenever the links have much more than two cells. The optimization models are solve using IBMIlogCplex on a laptop equipped with a 1.8 GHz processor and a 8.0 GB RAM. The optimization model (26) has 43 control variables and 1087 constraints and is solved in 0.06 s.

## IV. Robust Control for a Highway Network

The proposed framework is tested and evaluated at a network level. The study network is as shown in Figure 5. This network includes a 3-km section of highway CA-92 and a 3-km section of highway CA-101. The highways are divided into 8 links at the points where the number of lanes changes or where ramps connect the highway main lanes. The road characteristic parameters and traffic data are obtained from PeMS. Specifically, the length and number of lanes of each link are

$$L = [600, 600, 1200, 600, 600, 600, 1200, 600] \text{ (m)} \quad (27)$$
$$n_{lane} = [2, 2, 3, 3, 4, 5, 5, 5] \quad (28)$$

Each link is divided into 600-m segments. In order to have an unified segment length, the link lengths are approximated. For example, link 2 measured from Google map is 583 m which is less than a segment, so it is rounded up to 600 m. Such approximations are for simplicity, and in practical, the proposed model allows different segment sizes. In addition, the interchange between two highways is simplified as a junction within which the travel time is not considered. This assumption can also be relaxed by adding intermediate links between the junction and downstream links. The free flow speed is set to be 25 m/s according to the posted speed limits (55 mph). Based on the maximum capacity provided by the Vehicle Detector Stations (VDS), we let the capacity equal 1800 veh/h/lane. The jam density is 0.125 veh/m/lane, which falls within the common range of 0.115 to 0.155 veh/m/lane. Based on the triangular FD, the critical density and congestion speed are 0.02 veh/m/lane and $-4.76$ m/s, respectively.

For each link, the mean initial density is estimated in the same way as used in Section III-B. The locations of VDS for all links are shown in Figure 5. For the links that are not occupied by VDS, e.g., links 1 and 4, the closest VDS (400861 and 401343) are used. Flow and speed data between 9:00 am and 10:00 am in all weekdays of May 2019 are used. As a result, the mean initial density vector is expressed as

$$\bar{\rho} = [1.575, 0.975, 0.917, 0.352, 1.238, 0.730, 0.088, 0.084]$$
$$\odot \bar{\rho}_c, \quad (29)$$

where $\bar{\rho}_c$ is the critical density vector and $\odot$ denotes component-wise multiplication. We consider the uncertainties of links 3 and 7 so that the impact on both upstream and downstream links can be studied. It is assumed that the uncertainties follow normal distributions and the standard deviations are equal to $0.2\bar{\rho}(3)$ and $0.2\bar{\rho}(7)$ respectively.

### A. Network Modeling

We assume the traffic flow transition matrices at junctions are known, which is a commonly used assumption in traffic control models [16], [44]. Then, the traffic flow at a junction with $N_i$ incoming links and $N_o$ outgoing links can be modeled as,

$$\begin{bmatrix} q_{out} \\ q_{off} \end{bmatrix} = \begin{bmatrix} P^1 & P^2 \\ P^3 & 0 \end{bmatrix} \begin{bmatrix} q_{in} \\ q_{on} \end{bmatrix}, \quad (30)$$

where $q_{in}$ and $q_{out}$ are two column vectors denoting the incoming flows and outgoing flows at a junction node, respectively; $q_{on}$ and $q_{off}$ are two scalars representing the on-ramp and off-ramp flows; $P^1$ is a $N_o \times N_i$ matrix of which each element $P^1(i, j)$ means the proportion of the vehicles from incoming link $i$ going into link $j$; $P^2$ is a column vector with dimension of $N_o \times 1$ of which each element $P^2(i)$ means the proportion of the vehicles from on-ramp going into link $i$; $P^3$ is a row vector with dimension of $1 \times N_i$ of which each element $P^3(j)$ means the proportion of the vehicles from incoming link $j$ departing from the off-ramp. Usually, there is only one outgoing link at a node accepting on-ramp flows. Therefore, $P^2 = 1$ in most cases. In addition, we assumed no vehicles coming from an on-ramp would depart from the off-ramp at the same junction. The transition matrices are as follows,

$$P_2^1 = \begin{bmatrix} P(2,3) & P(6,3) \\ P(2,7) & P(6,7) \end{bmatrix} = \begin{bmatrix} 0.5 & 0.2 \\ 0.5 & 0.8 \end{bmatrix}$$
$$P_1^2 = P_3^2 = P_4^2 = P_5^2 = 1$$
$$P_3^3 = P_5^3 = 0.2, \quad (31)$$

The subscripts in (31) represent the junction nodes. Let $V$ and $L$ denote the sets of nodes and links, respectively; $q_{in}(i, j)$ and $q_{out}(i, j)$ are the inflows and outflows of link $j$ at time $i$. The traffic control problem for this highway network was modeled as,

$$\min_{\substack{q_{in}, q_{out}; \\ q_{on}, q_{off}}} -\sum_{i=1}^{n_{max}} \sum_{j=1}^{n_l} \Big( (q_{out}(i, j) + q_{in}(i, j))(n_{max} - i + 1) - \eta y(i) \Big)$$
$$s.t. \quad y(i) \geq n_{lane}(2) q_{out}(i, 6) - n_{lane}(6) q_{out}(i, 2), \quad \forall i$$
$$y(i) \geq n_{lane}(6) q_{out}(i, 2) - n_{lane}(2) q_{out}(i, 6), \quad \forall i$$

$$q_{on}(i, 1) \geq q_{out}(i, 1)/n_{lane}(1), \quad \forall i \in N$$
$$q_{on}(i, 2) \geq q_{out}(i, 3)/n_{lane}(3), \quad \forall i \in N$$
$$q_{on}(i, 3) \geq q_{out}(i, 5)/n_{lane}(5), \quad \forall i \in N$$
$$q_{on}(i, 4) \geq q_{out}(i, 7)/n_{lane}(7), \quad \forall i \in N$$
$$q_{out}(i, 4) \leq \psi'(\rho_4), \quad \forall i \in N$$
$$q_{out}(i, 8) \leq \psi'(\rho_8), \quad \forall i \in N$$
$$(38) - (42), \quad \forall j \in L$$
$$(30), \quad \forall v \in V$$
$$q_d^{out}(i, j) \geq 0, \quad q_d^{in}(i, j) \geq 0 \quad \forall i, j$$
$$q_{out}(i, j) \geq 0, \quad q_{in}(i, j) \geq 0 \quad \forall i, j \qquad (32)$$

The decision variables are the inflows and outflows of each link and the flows of on-ramps and off-ramps. The objective function consists of two terms. The first term is to maximize the sum of weighted outflows and inflows over all the links. A decreasing weight factor, $n_{max} - i + 1$, is added to force vehicles to move forward when space ahead is available. Otherwise, since maximizing unweighted inflows and outflows does not consider proceeding time, vehicles might have unnecessary stops that are not consistent with real-world conditions. The second term combined with the first two constraints is to add a penalty term if the outflows of two incoming links at a node are not proportional to their capacities. The parameter, $\eta$, is used to evaluate the trade-off between these two terms. A large $\eta$ imposes a low-level tolerance for the violation of proportionality while a small $\eta$ is inclined to allow the violation if it contributes to the first term more than the penalty. Let $\eta = 0.2$ in this section. The third to the sixth constraints set the outflow per lane on the merging links as the lower bound of the on-ramp inflows. Since the objective function only maximizes the flows on main-lane links, the on-ramp flow can be higher than this lower bound only when there is available space on the downstream link after it accommodates all vehicles from the upstream main-lane link. These constraints imply the priority of main-lane links is higher than on-ramps. The seventh and eighth constraints assume the traffic densities of downstream links of L4 and L8 are equal to their respective initial densities. Consequently, the upper bounds for their outflows are equal to the supply of those downstream links which can be expressed as $\psi'(\rho) = \rho_c v_f$ if $\rho \leq \rho_c$; $\psi'(\rho) = w(\rho - \rho_m)$ if $\rho > \rho_c$. As shown in Equation (29), the initial densities on links L4 and L8 are lower than their critical densities, so the upper bounds are their capacities in this example. Furthermore, we added the chance constraints (38)-(42) for each link and the transition model (30) for each node. Note that since the objective function maximizes the flows of every link in the network, the optimal solution implicitly avoids congestion to some extend so that the nodes can maintain a relatively high supply. For simplicity, instead of modeling the flows of the on-ramps through adding binary variables, we treat them as continuous variables. In addition, the simulation time is 500s and divided into 25 time steps, and the confidence level $1 - \alpha = 97.5\%$.

The optimization model (32) contains 550 control variables and 11036 constraints and is solved in 0.12 s. Fig. 6 shows the comparison of optimal flows of on-ramp 2 and on-ramp

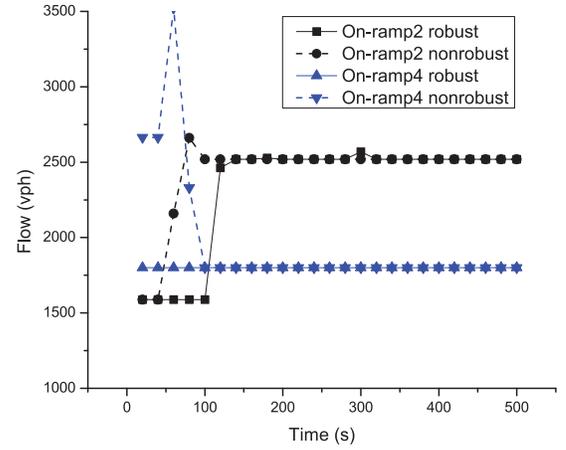

Fig. 6. Comparison of optimal flows of on-ramp 2 and on-ramp 4.

4 between the robust and non-robust models. As introduced at the beginning of this section, the test network is general, and the initial densities are estimated based on historical flow and speed data. Therefore, the density estimation does not show a regular and stable pattern. As a result, the simulation starts with shockwaves and fluctuate flows until the network reaches equilibrium, as shown in Figure 6. When considering the uncertainties on L3 and L7, to ensure the downstream nodes 3 and 5 can accommodate the merging flows with a high probability, the robust control model leads to lower on-ramp flows than the non-robust control model. Otherwise, the vehicles on both links might be blocked and generates a queue moving upstream, which results in a negative impact on the objective value. After a certain period, the influence of initial densities vanishes and the robust solution rises to the same value as the non-robust model.

Although only the control in a few of beginning steps is affected by the uncertainties, the control might impact the traffic state evolution in this network through the entire simulation. To validate the proposed model, let us consider a scenario in which the real initial densities on L3 and L7 exceed their mean values by one standard deviation. Then, forward simulations are conducted after implementing both the robust and non-robust optimal solutions on this scenario, and the results are compared in Section IV-C and Section IV-D.

### B. Sensitivity Analysis

The proposed approach requires the density in a spatial segment can be approximated as a constant, so this hyper parameter should be selected according to this standard. However, there is not such a rule for the time step selection. This section presents the influence of time step size on the control results. Let us use the robust optimization model (32) as an example. Since the optimal flow of on-ramp4 does not change during the simulation, we focus on the results of on-ramp 2. Figure 7 shows the optimal flow of on-ramp 2 with different step sizes. $n_{max}$ is the number of time steps.

It shows that our approach does not require time step to satisfy the CFL condition, which results in $n_{max} = 25$. Although a large time step may lose detail in the transition period, marked by ovals in Figure 7, the optimal solutions coincide with each other for most part of the simulation. Generally, a rule of thumb for the time step selection is the CFL condition, i.e. the time needed to cross the shortest cell at free flow speed. In addition, the computation time for these five scenarios is 0.02 s, 0.09 s, 0.11 s, 0.68 s, 6.735 s. Although the computation time is relative long when $n_{max} = 100$, the time step does not need to be so small in reality. In addition, based on the compatibility conditions (17)-(19), the number of constraints can be approximated as $n_l(2k_{max} \times n_{max} + n_{max}^2)$ where $n_l$ is the number of links. When $n_{max} = 25$, there are 11,036 constraints; when $n_{max} = 100$, there are 164,070 constraints. The computational burden from the number of temporal steps is much higher than the spatial segments. Therefore, our model is applicable for a larger network with a reasonable selection of time step and simulation time.

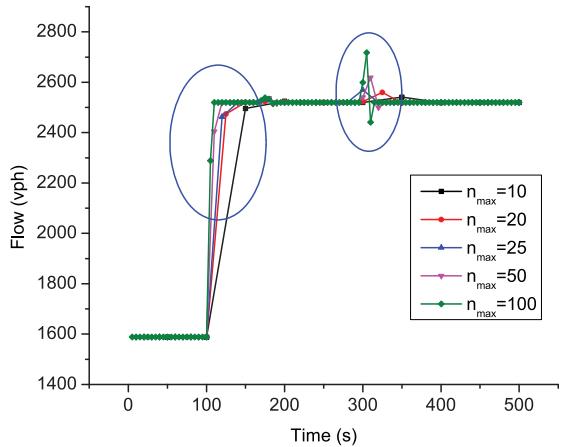

Fig. 7. Optimal flows of on-ramp 2 with different time steps.

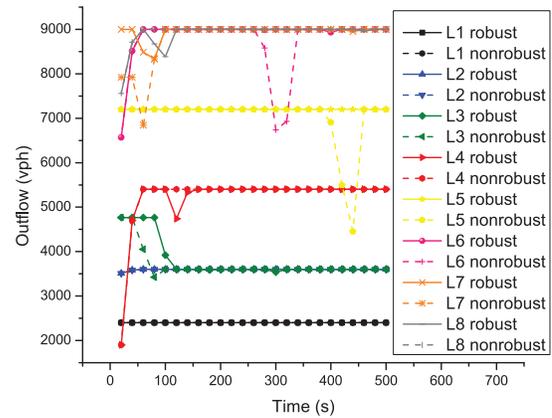

Fig. 8. Comparison of optimal outflow of link 2.

### C. On-Ramp Flow Control

In this first problem, we assume that we only control the on-ramp flows in the network. Since the initial densities on L1 and L4 are higher than their critical densities, the demand at both entry nodes are assumed to be equal to their capacities. The first term in the objective function, which maximizes the flows of every link in the network, requires the nodes to maintain a relatively high supply and implicitly avoids congestion to some extend. In case of an underestimate on initial densities, the non-robust optimal control allows excessive inflows from on-ramps to cause congestion and reduce outflows of links. Therefore, the comparison of link outflow is employed. After implementing the on-ramp control, the comparison of outflow on all links is shown in Figure 8. For each link, a solid line and a dashed line with the same color are used to represent the robust and non-robust results, respectively.

Following phenomena are observed in Figure 8. First, the outflow of boundary links L4 and L8 resulted from the robust control model is lower than the non-robust control model. The robust on-ramp control is derived from the model in which the confidence level is 97.5%, which means, in some sense, this control is optimal if the realization of density is approximately 1.96 standard deviation higher than the mean value. Therefore, the robust control is a little too conservative for this case, in which the initial densities are higher than the mean values by 1 standard deviation, and leads to undersaturated inflows at the beginning time steps on L4 and L8. Nevertheless, the chance constraints aim to make the optimal solution feasible with an expected probability rather than generate a solution that is optimal for all possible realizations. The second phenomenon is that, except for these two links, the robust control generates outflows higher or equal to the non-robust control on all other links. Although the average flow of on-ramp 2 and on-ramp 4 decreases by 103 vph and 159 vph, the overall increase of the outflow on these main lane links after considering the deduction on L4 and L8 is 550 vph. Therefore, the robust control solution is preferable in terms of gaining a better objective value. Third, the non-robust control generates a shockwave moving backward and causes a drop on outflows on L5 and L6 at the middle of simulation. Overall, this example shows that the robust control offers a higher throughput across the network.

### D. Boundary Network Control (Including On-Ramps)

In this problem, we implement the optimal inflows of L1 and L5, in addition to the on-ramps, on the scenario. The evolution of densities across all links is shown in Fig. 9. The axes denote spatial coordinates, and the unit is m. Since we fixed the on-ramp flows in this case, the flows from the on-ramps have a higher priority than the highway links to be served. Since the initial densities are underestimated in the non-robust control method, it sends excessive vehicles from on-ramps. Consequently, L4 and L8 do not have enough space to serve the sending flows from L3 and L7. It clearly shows the generation of congestion on L3 and L7, and its evolution across the network resulted from the non-robust control implementation. In contrast, no such congestion is induced from the robust control. L1 is highly congested at the beginning, and the objective function forces as many vehicles to enter as it can provide space for. Therefore, it is congested during the entire simulation for both cases. This case study depicts that the robust control is able to avoid congestion when the initial density is underestimated.

## V. MONTE CARLO SIMULATION WITHOUT RELAXATION

### A. Difficulty Without Relaxation

In the stochastic optimization model derived in Appendix A, some constraints were relaxed, since only $\rho(k)$ for the constraints involving $M_{M_k}$ was considered as a random variable while all of other $\rho(i)'s, i \in \{1, 2, \ldots, k-1\}$ were regarded as fixed values with their corresponding means. In this section, the difficulty of converting the chance constraints into linear form without relaxation was explained. After this, the Monte Carlo simulation was executed to remove the relaxation. At the end of this section, the optimal objective values (total outflows) from the Monte Carlo simulation were compared to the results from the last section. For simplicity, the Moskowitz solutions from the initial conditions (11) can be expressed as:

$$M_{M_k}(t, x) = \begin{cases} f_1(\rho(i)), & \text{if } \rho(k) \geq \rho_c, i = 1, 2, \ldots, K_{max} \\ f_2(\rho(i)), & \text{if } \rho(k) < \rho_c, i = 1, 2, \ldots, K_{max} \end{cases}$$
(33)

where $f_1$ and $f_2$ indicate two linear functions and $K_{max}$ is the number of the initial condition segments. As a result, the typical chance constraint involving $M_{M_k}$ can be expressed as,

$$Pr(M_{M_k}(t, x) \geq g(q))$$
$$= Pr(f_1(\rho(i)) \geq g(q), \rho(k) \geq \rho_c)$$
$$+ Pr(f_2(\rho(i)) \geq g(q), \rho(k) < \rho_c) \quad (34)$$

where $g$ is a linear function of boundary conditions.

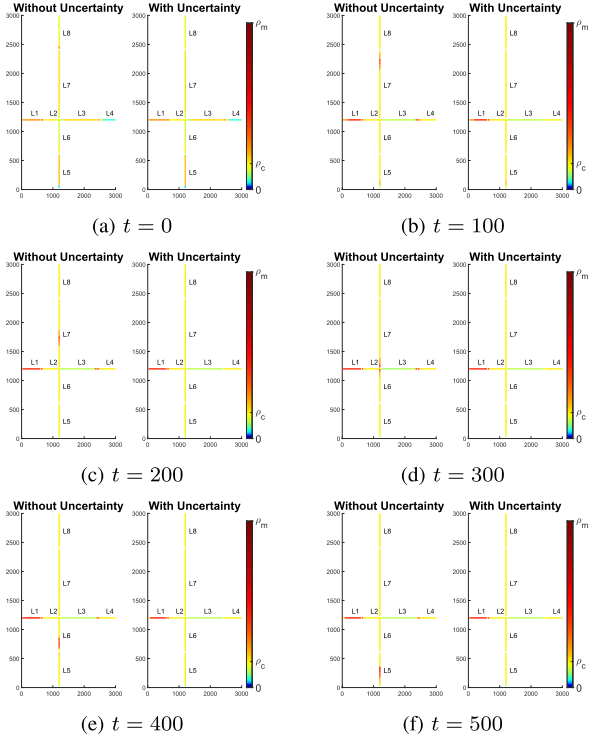

Fig. 9. Evolution of the density across the network, with classical control (left) and robust control (right).

If all of the initial conditions are independently normally distributed, $(f(\rho(i)), \rho(k))$ is subject to a bivariate normal distribution $n(\boldsymbol{\mu}, \boldsymbol{\Sigma})$,

$$\boldsymbol{\mu} = [\mu_{f(\rho(i))}, \mu_{\rho(k)}],$$
$$\boldsymbol{\Sigma} = \begin{bmatrix} Var(f(\rho(i))) & Cov(f(\rho(i), \rho(k))) \\ Cov(f(\rho(i), \rho(k))) & Var(\rho(k)) \end{bmatrix} \quad (35)$$

Although the pdf of a bivariate normal distribution can be obtained, there is no closed-form of the corresponding cumulative distribution function (cdf). Therefore, for a bivariate normal distribution, the chance constraints cannot be expressed as linear constraints by the same method.

### B. Monte Carlo Simulation

To validate our relaxed model, Monte Carlo simulations were used to convert the chance constraints into a linear form. The algorithm for constraints $Pr(M_{M_k}(t, x) \geq \gamma(t, x)) \geq 1 - \alpha$ is as follows.

*Step 1:* Generate $N$ random numbers from the normal distribution for each initial condition segment. In this paper, $N = 1000$ and $\rho(k_i)$ is the $i$th number for the $k$th segment.

*Step 2:* Calculate $M_{M_k(i)}(t, x), i = 1, 2, \ldots N$ using the $i$th number from each segment from Step 1.

*Step 3:* Sort $M_{M_k(i)}(t, x)$ into ascending order. Find the corresponding critical value. For example, if the confidence level is 97.5%, then the critical value should be $M_{M_k(25)}(t, x)$ in the ordered sequence.

*Step 4:* Replace the constraints involving $M_{M_k}$ and $\gamma$ of $Pr(M_{M_k}(t, x) \geq \gamma(t, x)) \geq 1 - \alpha$ with $M_{M_k(N\alpha)}(t, x) \geq \gamma(t, x)$.

Because the downstream condition of $\beta_n(t, x)$ is a function of the initial conditions, the algorithm for the constraints $Pr(M_{M_k}(t, x) \geq \beta(t, x)) \geq 1 - \alpha$ is as follows.

*Step 1:* Generate $N$ random numbers from the normal distribution for each initial condition segment. $\rho(k_i)$ is the $i$th number for the $k$th segment.

*Step 2:* Calculate $M_{M_k(i)}(t, x) + \sum_{k=1}^{k=k_{max}} \rho(k_i)X, i = 1, 2, \ldots N$ using the $i$th number from each segment from Step 1.

*Step 3:* Sort $M_{M_k(i)}(t, x) + \sum_{k=1}^{k=k_{max}} \rho(k_i)X$ in an ascending order. Use $(M_{M_k}(t, x) + \sum_{k=1}^{k=k_{max}} \rho(k)X)_{N\alpha}$ to represent the $N\alpha$th element in the ordered sequence.

*Step 4:* Replace the constraints involving $M_{M_k}$ and $\beta$ of $Pr(M_{M_k}(t, x) \geq \beta(t, x)) \geq 1 - \alpha$ with $(M_{M_k}(t, x) + \sum_{k=1}^{k=k_{max}} \rho(k)X)_{N\alpha} \geq \beta(t, x) + \sum_{k=1}^{k=k_{max}} \rho(k)X$.

There is no need to modify any constraints not involving $M_{M_k}(t, x)$ because there is no relaxation in those constraints. Since there is no constraint for the capacity downstream, the time period in which the outflow is impacted by the initial density could be approximated by the link length divided by the free flow speed, which is 7 steps in this example. The comparison of averaged outflow over this impacted time period for the problem (24) is shown in Fig. 10.

Since the objective function is to maximize the total outflows, the optimal value from Monte Carlo simulation is smaller than the relaxed LP. The error increases with the confidence level and the standard deviation, and the largest error percent is 15%. Therefore, our relaxed stochastic program fits well with the Monte Carlo simulation. To explain the reason why the optimal solution from the relaxed LP has such a high accuracy, let us delve into some constraints. The errors that originated from the relaxed constraints are mitigated by these accurate constraints. (39) shows that, compared to the Monte Carlo simulation, the relaxed constraints generate a larger upper bound on the total outflows, and this could lead to an error approximating $\frac{z_{1-\alpha}\sqrt{\sum_{i=k+1}^{k_{max}}\sigma(i)}}{\sum_{i=k+1}^{k_{max}}\rho(i)+z_{1-\alpha}\sqrt{\sum_{i=k+1}^{k_{max}}\sigma^2(i)}}$ in percentage. Although the constraints (38)-(40) are relaxed, the constraints (41)-(42) are accurate. For example, when $n=1$, $p=7$, and the first constraint in (41) could be written as

$$-q_{in}(1)(0.57T) + \sum_{1}^{7} q_{out}(i)T \leq \left(\sum_{k=1}^{k=k_{max}} \rho_k - z_{1-\alpha}\sqrt{\sum_{k=1}^{k=k_{max}} \sigma_k^2}\right) \quad (36)$$

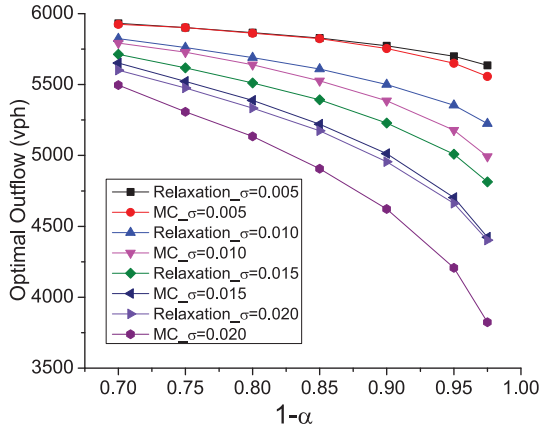

Fig. 10. Comparison of optimally total outflows.

This constraint is accurate and produces a stricter upper bound. Therefore, the couple effect of these constraints ensures the accuracy of our model.

## VI. CONCLUSION

This paper developed a robust traffic control model to consider the uncertainties in initial densities by using chance constraints. The proposed model is relaxed and linearized to be solved efficiently. Following this, case studies for both a single highway link and a network show the benefits of the proposed model. In addition, Monte Carlo simulation was used to verify the accuracy of the relaxed model.

A set of individual chance constraints was used for the proposed model. In future work, considering joint chance constraints is a promising topic. In addition, finding a stochastic programming formulation to deal with the uncertainty in model parameters will be another interesting research direction.

## APPENDIX

### A. Linear Expression of the Chance Constraints

The following procedure shows the derivation of the deterministic version for the chance constraint $Pr\{M_{M_k}(pT,\xi) \geq \gamma_p(pT,\xi)\} \geq 1-\alpha$ in detail, and the integrated deterministic constraints were obtained using the same method.

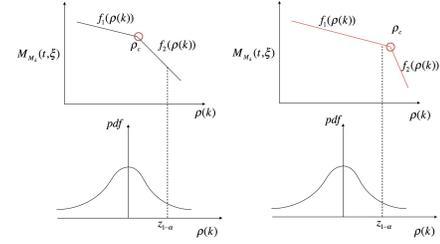

Fig. 11. Solution at upstream from initial condition and associated initial density distribution.

From (11), the Moskowitz solution upstream from the initial condition (6) can be explicitly expressed as:

$$M_{M_k}(t,\xi) = \begin{cases} +\infty, & \text{if } t \leq \frac{\xi-(k-1)X}{w} \\ -\sum_{i=1}^{k-1}\rho(i)X + \rho_c(tv_f + (k-1)X - \xi), \\ & \text{if } t \geq \frac{\xi-(k-1)X}{w} \text{ and } \rho(k) \leq \rho_c \\ -\sum_{i=1}^{k-1}\rho(i)X + \rho(k)(tw+(k-1)X-\xi) \\ -\rho_m tw, \\ & \text{if } \frac{\xi-(k-1)X}{w} \leq t \leq \frac{\xi-kX}{w} \text{ and } \rho(k) \geq \rho_c \\ -\sum_{i=1}^{k}\rho(i)X + \rho_c(tw+kX-\xi) - \rho_m tw, \\ & \text{if } t \geq \frac{\xi-kX}{w} \text{ and } \rho(k) \geq \rho_c \end{cases} \quad (37)$$

From (37), it is known that $M_{M_k}(t,\xi)$ is a nonincreasing and piecewise linear function of $\rho(k)$, as shown in Fig. 11. Then the corresponding chance constraint was simply divided into two situations:

(i). $\rho_c \leq \rho_k + z_{1-\alpha}\sigma_k$ as shown left in Fig. 11. We should convert the chance constraint to $f_2(\rho_k + z_{1-\alpha}\sigma_k) \geq \gamma_p(pT,\xi)$ $\forall k \in K$, $\forall p \in N$;

(ii). $\rho_c \geq \rho_k + z_{1-\alpha}\sigma_k$ as shown right in Fig. 11. We should convert the chance constraint to $f_1(\rho_k + z_{1-\alpha}\sigma_k) \geq \gamma_p(pT,\xi)$ $\forall k \in K$, $\forall p \in N$.

where $z_{1-\alpha}$ is defined as z score such that $P(\rho(k) \leq \rho_k + z_{1-\alpha}\sigma_k) = 1-\alpha$.

Substituting the expressions of $M_{M_k}(t,\xi)$ (37) and $\gamma_p(pT,\xi)$ (7) into the inequalities above leads to the following linear and deterministic constraint:

$$\begin{cases} -\sum_{i=1}^{k-1}\rho(i)x + \rho_c(pTv_f + (k-1)x - \xi) \\ \geq \sum_{i=1}^{p} q_{in}(i)T, \\ \quad \text{if } t \geq \frac{\xi-(k-1)X}{w} \\ \quad \text{and } \rho_k + z_{1-\alpha}\sigma_k \leq \rho_c \\ -\sum_{i=1}^{k-1}\rho(i)X + (\rho(k)+z_{1-\alpha}\sigma_k)(tw+(k-1)X-\xi) \\ -\rho_m tw \geq \sum_{i=1}^{p} q_{in}(i)T, \\ \quad \text{if } \frac{\xi-(k-1)X}{w} \leq t \leq \frac{\xi-kX}{w}, \\ \quad \text{and } \rho_k + z_{1-\alpha}\sigma_k \geq \rho_c \\ -\sum_{i=1}^{k-1}\rho(i)X - (\rho(k)+z_{1-\alpha}\sigma_k)X + \rho_c(tw+kX-\xi) \\ -\rho_m tw \geq \sum_{i=1}^{p} q_{in}(i)T, \\ \quad \text{if } t \geq \frac{\xi-kX}{w}, \\ \quad \text{and } \rho_k + z_{1-\alpha}\sigma_k \geq \rho_c \end{cases} \quad (38)$$

For simplicity, only $\rho(k)$ for the constraints involving $M_{M_k}$ was considered as a random variable, all other $\rho(i)'s, i \in \{1,2,\ldots,k-1\}$ were still regarded as fixed values with their corresponding means. In Section V, the complexity of regarding all of the $\rho(i)'s, i \in \{1, 2, \ldots, k-1\}$ as random variables will be explained.

The other constraints were found in this same manner. The deterministic expression for the rest of the constraints in (17) is shown in (39) and (40), and the deterministic expression of (18) and (19) are shown in (41) and (42), respectively,

$$\begin{cases} \sum_{i=k+1}^{k_{max}} \rho(i)X + (\rho_k - z_{1-\alpha}\sigma_k)(pTv_f + kX - \chi) \\ \geq \sum_{i=1}^{p} q_{out}(i)T, \\ \qquad \text{if } \frac{\chi-kX}{v_f} \leq pT \leq \frac{\chi-(k-1)X}{v_f}, \\ \qquad \text{and } \rho_k - z_{1-\alpha}\sigma_k \leq \rho_c \\ \sum_{i=k+1}^{k_{max}} \rho(i)X + (\rho_k - z_{1-\alpha}\sigma_k)X + \\ \rho_c(pTv_f + (k-1)X - \chi) \\ \geq \sum_{i=1}^{p} q_{out}(i)T, \\ \qquad \text{if } pT \geq \frac{\chi-(k-1)X}{v_f}, \\ \qquad \text{and } \rho_k - z_{1-\alpha}\sigma_k \leq \rho_c \\ \sum_{i=k+1}^{k_{max}} \rho(i)X + \rho_c(pTv_f + kX - \chi) - \rho_m pTw \\ \geq \sum_{i=1}^{p} q_{out}(i)T, \\ \qquad \text{if } pT \geq \frac{\chi-kX}{v_f}, \\ \qquad \text{and } \rho_k - z_{1-\alpha}\sigma_k \geq \rho_c \end{cases} \tag{39}$$

$$\begin{cases} -\sum_{i=1}^{k-1}\rho_i x + \rho_c(\frac{\xi-(k-1)X}{w}v_f + (k-1)x - \xi) \\ \geq \sum_{i=1}^{p-1} q_{in}(i)T + q_{in}(p)(\frac{\xi-(k-1)X}{w} - (p-1)T), \\ \qquad \text{if } (p-1)T \leq \frac{\xi-(k-1)X}{w} \leq pT, \\ \qquad \text{and } \rho_k + z_{1-\alpha}\sigma_k \leq \rho_c \\ -\sum_{i=1}^{k-1}\rho_i x - \rho_m(\xi-(k-1)X) \\ \geq \sum_{i=1}^{p-1} q_{in}(i)T + q_{in}(p)(\frac{\xi-(k-1)X}{w} - (p-1)T), \\ \qquad \text{if } (p-1)T \leq \frac{\xi-(k-1)X}{w} \leq pT, \\ \qquad \text{and } \rho_k + z_{1-\alpha}\sigma_k \geq \rho_c \\ \sum_{i=k+1}^{k_{max}} \rho_i x \\ \geq \sum_{i=1}^{p-1} q_{out}(i)T + q_{out}(p)(\frac{\chi-kX}{v_f} - (p-1)T), \\ \qquad \text{if } (p-1)T \leq \frac{\chi-kX}{v_f} \leq pT, \\ \qquad \text{and } \rho_k + z_{1-\alpha}\sigma_k \leq \rho_c \\ \sum_{i=k+1}^{k_{max}} \rho_i x + \rho_c(\frac{\chi-kX}{v_f}w + kX - \chi) \\ \geq \sum_{i=1}^{p-1} q_{out}(i)T + q_{out}(p)(\frac{\chi-kX}{v_f} - (p-1)T), \\ \qquad \text{if } (p-1)T \leq \frac{\chi-kX}{v_f} \leq pT, \\ \qquad \text{and } \rho_k + z_{1-\alpha}\sigma_k \geq \rho_c \end{cases} \tag{40}$$

$$\begin{cases} -\sum_{i=1}^{n-1} q_{in}(i)T - q_{in}(n)(pT - \frac{\chi-\xi}{v_f} - (n-1)T) \\ +\sum_{i=1}^{p} q_{out}(i)T \leq (\sum_{k=1}^{k_{max}} \rho_k - z_{1-\alpha}\sqrt{\sum_{k=1}^{k_{max}} \sigma_k^2})X \\ \qquad \text{if } (n-1)T + \frac{\chi-\xi}{v_f} \leq pT \leq nT + \frac{\chi-\xi}{v_f} \\ -\sum_{i=1}^{n} q_{in}(i)T - \rho_c v_f(pT - \frac{\chi-\xi}{v_f} - nT) + \sum_{i=1}^{p} q_{out}(i)T \\ \leq (\sum_{k=1}^{k_{max}} \rho_k - z_{1-\alpha}\sqrt{\sum_{k=1}^{k_{max}} \sigma_k^2})X \\ \qquad \text{if } pT \geq nT + \frac{\chi-\xi}{v_f} \\ -\sum_{i=1}^{n} q_{in}(i)T + \sum_{i=1}^{p-1} q_{out}(i)T \\ +q_{out}(p)(nT + \frac{\chi-\xi}{v_f} - (p-1)T) \\ \leq (\sum_{k=1}^{k_{max}} \rho_k - z_{1-\alpha}\sqrt{\sum_{k=1}^{k_{max}} \sigma_k^2})X \\ \qquad \text{if } (p-1)T \leq nT + \frac{\chi-\xi}{v_f} \leq pT \end{cases} \tag{41}$$

$$\begin{cases} \sum_{i=1}^{n-1} q_{out}(i)T + q_{out}(n)(pT - \frac{\xi-\chi}{w} - (n-1)T) \\ -\rho_m(\xi-\chi) - \sum_{i=1}^{p} q_{in}(i)T \\ \geq (\sum_{k=1}^{k_{max}} \rho_k + z_{1-\alpha}\sqrt{\sum_{k=1}^{k_{max}} \sigma_k^2})X \\ \qquad \text{if } (n-1)T + \frac{\xi-\chi}{w} \leq pT \leq nT + \frac{\xi-\chi}{w} \\ \sum_{i=1}^{n} q_{out}(i)T + \rho_c v_f(pT - \frac{\xi-\chi}{w} - nT) \\ -\sum_{i=1}^{p} q_{in}(i)T \\ \geq (\sum_{k=1}^{k_{max}} \rho_k + z_{1-\alpha}\sqrt{\sum_{k=1}^{k_{max}} \sigma_k^2})X \\ \qquad \text{if } pT \geq nT + \frac{\xi-\chi}{w} \\ \sum_{i=1}^{n} q_{out}(i)T + \rho_c v_f(\frac{\xi-\chi}{w} - \frac{\xi-\chi}{v_f}) \\ -\sum_{i=1}^{p-1} q_{in}(i)T + q_{in}(p)(nT + \frac{\xi-\chi}{w} - (p-1)T) \\ \geq (\sum_{k=1}^{k_{max}} \rho_k + z_{1-\alpha}\sqrt{\sum_{k=1}^{k_{max}} \sigma_k^2})X \\ \qquad \text{if } (p-1)T \leq nT + \frac{\xi-\chi}{w} \leq pT \end{cases} \tag{42}$$

### B. Proof of the Relationship Between the Total Outflow and the Differentials of Outflows Across Time Steps

Assume $q'_{out}(i)$ is the optimal solution of the objective function $max \sum_{i=1}^{n_{max}} q_{out}(i)$.

For any feasible solution of (24), the objective function is

$$\sum_{i \in Q^{++}} (h+2)q_{out}(i) + \sum_{i \in Q^+} (h+1)q_{out}(i) + \sum_{i \in Q} hq_{out}(i)$$
$$+ \sum_{i \in Q^-} (h-1)q_{out}(i) + \sum_{i \in Q^{--}} (h-2)q_{out}(i) \tag{43}$$

where

$$Q^{++} = \{i : q_{out}(i) < q_{out}(i+1) \cap q_{out}(i) < q_{out}(i-1)\} \tag{44}$$

$$Q^+ = \{i : q_{out}(i) < q_{out}(i+1) \cap q_{out}(i) = q_{out}(i-1) \\ \text{or } q_{out}(i) = q_{out}(i+1) \cap q_{out}(i) < q_{out}(i-1)\} \tag{45}$$

$$Q^{--} = \{i : q_{out}(i) > q_{out}(i+1) \cap q_{out}(i) > q_{out}(i-1)\} \tag{46}$$

$$Q^- = \{i : q_{out}(i) > q_{out}(i+1) \cap q_{out}(i) = q_{out}(i-1) \\ \text{or } q_{out}(i) = q_{out}(i+1) \cap q_{out}(i) > q_{out}(i-1)\} \tag{47}$$

$$Q = \{i : i \in I - Q^{++} - Q^+ - Q^{--} - Q^-\} \tag{48}$$

$q_{out}(1)$ and $q_{out}(n_{max})$ are in $Q^+$, $Q$ or $Q^-$ based on the same logic. When $h > 2$, the coefficient for any time step $i$ is positive.

Assume the optimal solution of (43) is $q^*_{out}(i)$ and the total outflow $\sum_{i=1}^{n_{max}} q^*_{out}(i)$ is less than $\sum_{i=1}^{n_{max}} q'_{out}(i)$, there must exist at least one time step such that $q^*_{out}(i) < q'_{out}(i)$.

Since the coefficients are all positive, let $q'_{out}(i)$ replace $q^*_{out}(i)$, the new objective value of (43) will be larger than $\sum_{i=1}^{n_{max}} q^*_{out}(i)$. Therefore, $q^*_{out}(i)$ is not the optimal solution. Therefore, the optimally total outflow of $\sum_{i=1}^{n_{max}} q^*_{out}(i)$ equals $\sum_{i=1}^{n_{max}} q'_{out}(i)$.

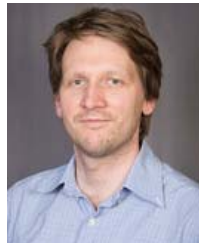

**Christian Claudel** received the M.S. degree in plasma physics from the École normale supérieure de Lyon in 2004, and the Ph.D. degree in electrical engineering and computer sciences from the Department of Electrical Engineering and Computer Sciences, UC Berkeley, in 2010. He was an Assistant Professor of electrical engineering with the King Abdullah University of Science and Technology. He is currently an Associate Professor of civil architectural and environmental engineering with The University of Texas at Austin. His research interests include the control and estimation of distributed parameter systems, wireless sensor networks, and unmanned air vehicles.

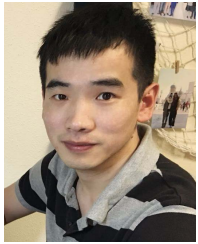

**Hao Liu** received the M.S. degree in statistics and the Ph.D. degree in civil, architectural, and environmental engineering from The University of Texas at Austin in 2018 and 2020, respectively. His research interests include traffic flow modeling and optimization.

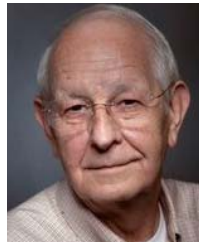

**Randy Machemehl** is currently a Professor of transportation engineering with The University of Texas at Austin and the former Director of the Center for Transportation Research. He was in private engineering practice, as a Staff Member of the Wilbur Smith and Associates, before joining the faculty of The University of Texas at Austin in 1978. His research interests include transportation system operations, public transportation systems planning and design, traffic data acquisition, traffic simulation, optimization, and bicycle safety.